\documentclass[12pt,reqno]{amsart}
\usepackage{amssymb}
\usepackage[all]{xy}

\newtheorem{thm}{Theorem}[section]
\newtheorem{prop}[thm]{Proposition}
\newtheorem{lem}[thm]{Lemma}
\newtheorem{cor}[thm]{Corollary}

\theoremstyle{definition}
\newtheorem{defn}[thm]{Definition}
\newtheorem{remark}[thm]{Remark}

\newtheorem{example}[thm]{Example}

\newtheorem{reduction}[thm]{Reduction}
\newtheorem*{ack}{Acknowledgment}

\numberwithin{equation}{section}

\newcommand{\mo}{\mathopen<}
\newcommand{\mc}{\mathclose>}
\newcommand{\pfistero}{\ll}
\newcommand{\pfisterc}{\gg}
\newcommand{\ed}{\operatorname{ed}}
\newcommand{\Alg}{\operatorname{Alg}}

\newcommand{\Sym}{{\operatorname{\mathcal{S}}}}
\newcommand{\rank}{\operatorname{rank}}

\newcommand{\im}{\operatorname{Im}}

\newcommand{\Char}{\operatorname{char}}
\newcommand{\Ker}{\operatorname{Ker}}

\newcommand{\Hom}{\operatorname{Hom}}

\newcommand{\U}{\operatorname{U}}
\newcommand{\UD}{\operatorname{UD}}

\newcommand{\Ext}{\operatorname{Ext}}
\newcommand{\Gal}{\operatorname{Gal}}

\newcommand{\pr}{\operatorname{pr}}
\newcommand{\trdeg}{\operatorname{trdeg}}
\newcommand{\tr}{\operatorname{Tr}}

\newcommand{\onto}{\twoheadrightarrow}
\newcommand{\mono}{\rightarrowtail}
\newcommand{\into}{\hookrightarrow}
\newcommand{\tto}{\longrightarrow}

\newcommand{\Mat}{{\operatorname{M}}}
\newcommand{\Mn}{\Mat_n}
\newcommand{\PGLn}{{\operatorname{PGL}}_n}

\newcommand{\wed}[1]{\textstyle\bigwedge^2{#1}}
\newcommand{\sym}[1]{\operatorname{\mathsf{Sym}}^2{#1}}
\newcommand{\ind}[2]{\!\!\uparrow_{#1}^{#2}}

\newcommand{\cG}{\mathcal{G}}
\newcommand{\cH}{\mathcal{H}}
\newcommand{\cN}{\mathcal{N}}
\newcommand{\cR}{\mathcal{R}}
\newcommand{\cF}{\mathcal{F}}

\newcommand{\bt}{\mathbf{t}}

\newcommand{\g}{\overline{g}}

\newcommand{\x}{\overline{g_1}}
\newcommand{\y}{\overline{g_2}}

\newcommand{\oGH}{\omega(\GH)}

\newcommand{\GH}{\cG/\cH}

\newcommand{\bbZ}{\mathbb{Z}}

\newcommand{\bbQ}{\mathbb{Q}}

\newcommand{\BB}{\mathbf{B}}
\newcommand{\Br}{\mathbf{B}}

 \newdir{ _{(}}{{}*!/-3pt/\dir_{(}}
 \newdir{ ^{(}}{{}*!/-3pt/\dir^{(}}

\begin{document}

\title[Fields of definition]{Fields of definition for division algebras}

\author{M. Lorenz}
\address{Department of Mathematics, Temple University,
 Philadelphia, PA 19122-6094}
\email{lorenz@math.temple.edu}
\urladdr{http://www.math.temple.edu/$\stackrel{\sim}{\phantom{.}}$lorenz}
\thanks{M. Lorenz was supported in part by NSF grant
DMS-9988756}

\author{Z. Reichstein}
\address{Department of Mathematics, University of British Columbia,
 Vancouver, BC, Canada, V6T 1Z2}
\email{reichst@math.ubc.ca}
\urladdr{http://www.math.ubc.ca/$\stackrel{\sim}{\phantom{.}}$reichst}
\thanks{Z. Reichstein was supported in part by NSF grant DMS-9801675 and
an NSERC research grant}

\author{L. H. Rowen}
\address{Department of Mathematics and Computer Science,
Bar-Ilan University, \newline
Ramat-Gan  52900, Israel}
\email{rowen@macs.biu.ac.il}
\thanks{L. H. Rowen was supported by the Israel Science Foundation,
founded by the Israel Academy
of Sciences and Humanities - Center of Excellence Program no. ~8007/99-3}

\author{D. J. Saltman}
\address{Department of Mathematics, University of Texas,
 Austin, TX 78712}
\email{saltman@math.utexas.edu}
\thanks{D. J. Saltman was supported in part by NSF grant DMS-9970213}

\begin{abstract}
Let $A$ be a finite-dimensional division algebra
containing a base field $k$ in its center $F$.
We say that $A$ is defined over a subfield $F_0$
if there exists an $F_0$-algebra $A_0$ such that
$A = A_0 \otimes_{F_0} F$.
We show that: (i) In many cases $A$ can be defined over
a rational extension of $k$. (ii) If $A$ has odd degree $n \geq 5$,
then $A$ is defined over a field $F_0$ of transcendence
degree $\leq {\frac{1}{2}}(n-1)(n-2)$ over $k$.
(iii) If $A$ is a $\bbZ/m \times \bbZ/2$-crossed product for some $m \geq 2$
(and in particular, if $A$ is any algebra of degree $4$)
then $A$ is Brauer equivalent to a tensor product of two symbol algebras.
Consequently, $\Mat_m(A)$ can be defined
over a field $F_0$ such that $\trdeg_k(F_0) \leq 4$.
(iv) If $A$ has degree 4 then the trace form of $A$ can
be defined over a field $F_0$ of transcendence degree $\leq 4$.
(In (i), (iii) and (iv) we assume that the center of $A$ contains
certain roots of unity.)
\end{abstract}

\maketitle
\tableofcontents

\section{Introduction}

Let $F$ be a field and $A$ a finite-dimensional $F$-algebra.
We say that $A$ is
defined over a subfield $F_0 \subset F$ if there exists an
$F_0$-algebra $A_0$ such that $A = A_0 \otimes_{F_0} F$.

Throughout this paper we will assume that $A$ is a finite-dimensional
central simple $F$-algebra, and $F$ (and $F_0$) contain a base field $k$.

\subsection{Parameter reduction}
If $A$ is defined over $F_0$ and $\trdeg_k(F_0) < \trdeg_k(F)$ then
the passage from $A$ to $A_0$ may be viewed as ``parameter reduction" in
$A$.
This leads to the following natural question:

\begin{quote}
What is the smallest value of $\trdeg_k(F_0)$ such that $A$ is defined
over $F_0$?
\end{quote}

\noindent
This number is clearly finite; we shall denote it by $\tau(A)$.
Of particular interest to us will be the case where $A = \UD(n)$ is
the universal division algebra of degree $n$ and $F = Z(n)$ is the center
of $\UD(n)$. Recall that $\UD(n)$ is
the subalgebra of $\Mn(k(x_{ij}, y_{ij}))$ generated,
as a $k$-division algebra, by two generic $n \times n$-matrices
$X = (x_{ij})$ and $Y = (y_{ij})$,
where $x_{ij}$ and $y_{ij}$ are $2n^2$ independent commuting variables over
$k$;
see, e.g.,~\cite[Section II.1]{procesi}, \cite[Section 3.2]{rowen-pi}
or~\cite[Section 14]{saltman}.
We shall write $d(n)$ for $\tau(\UD(n))$.
It is easy to show that $d(n) \geq \tau(A)$ for any central
simple algebra $A$ of degree $n$ whose center contains $k$; see
Remark~\ref{rem-ud(n)} (cf. also~\cite[Lemma 9.2]{reichstein1}).
In other words, every central simple algebra of degree $n$
can be ``reduced to at most $d(n)$ parameters".

To the best of our knowledge, the earliest attempt
to determine the value of $d(n)$ is due to
Procesi, who showed that $d(n) \leq n^2$; see~\cite[Thm. 2.1]{procesi}.
If $n=2,3$ or $6$ and $k$ contains a primitive
$n^{\text{th}}$ root of unity then $d(n) = 2$,
because $\UD(n)$ is cyclic for these $n$ and we can take $A_0$ to be a
symbol algebra;
cf.~\cite[Lemma 9.4]{reichstein1}. Rost~\cite{rost} recently proved
that
\begin{equation} \label{e.rost}
d(4) = 5 \, .
\end{equation}
For other $n$ the exact value of $d(n)$ is not known.
However, the following inequalities hold:
\begin{gather}
  d(n)  \leq  n^2 -3n + 1 \quad \text{\rm if $n \geq 4$}
 \quad \text{\rm \cite{lemire}} \, ,  \nonumber \\
 d(n)  \leq  d(nm) \leq d(n) + d(m) \quad
 \text{\rm if $(n, m) = 1$} \quad
 \text{\rm \cite[Sect.~9.4]{reichstein1}} \, , \nonumber \\
 d(n^r)  \geq   2r \quad
 \text{\rm \cite[Theorem 16.1]{reichstein2}} \, , \\
 d(n)  \leq  \tfrac{1}{2}(n-1)(n-2) + n  \quad \text{\rm for odd $n$}
 \quad
 \text{\rm \cite{rowen-brauer}; cf.~\cite[Sect.~9.3]{reichstein1}} \, .
\nonumber
\end{gather}
In this paper we will sharpen the last inequality
by showing that $$d(n) \leq \frac{1}{2}(n-1)(n-2)$$
for every odd $n \geq 5$. Moreover, in $\UD(n)$, reduction to this number
of parameters can be arranged in a particularly nice fashion:

\begin{thm} \label{thm1} Let $n \geq 5$ be an odd integer, $\UD(n)$
the universal division algebra of degree $n$, and
$Z(n)$ its center.  Then there exists a subfield $F$ of $Z(n)$ and
a division algebra $D$ of degree $n$ with center $F$ such that
\begin{enumerate}
\item
$\UD(n) = D \otimes_F Z(n)$,
\item
$\trdeg_k(F) = \frac{1}{2}(n-1)(n-2)$, and
\item
$Z(n)$ is a rational extension of $F$.
\end{enumerate}
In particular, $d(n) \le \frac{1}{2}(n-1)(n-2)$.
\end{thm}

We remark that in the language of~\cite{reichstein1}, the last assertion
of Theorem~\ref{thm1} can be written as
\begin{equation} \label{e.ed-reference}
\ed(\PGLn) \le \frac{1}{2}(n-1)(n-2).
\end{equation}

\subsection{Rational fields of definition}

Another natural question is whether or not a given central simple algebra
$A$ can be defined over a rational extension of $k$. We give the following
partial answer to this question.

\begin{thm} \label{thm.rat}
Let $A$ be a finite-dimensional
central simple algebra with center $F$ and let
$t_1, t_2, \dots$ be algebraically independent central indeterminates
over $F$.
\begin{enumerate}

\item
Assume $\deg(A) = 2^i p_1 \dots p_r$, where $i = 0$, $1$ or $2$ and
$p_1, \dots, p_r$ are distinct odd primes. Then for $s \gg 0$
the algebra $A(t_1, \dots, t_s)$ is defined over a rational
extension of $k$.

\item
Suppose the center of $A$ contains a primitive $e^{\text{th}}$ root of
unity,
where $e$ is the exponent of $A$. Then there exists an $r \geq 1$ such
that for $s \gg 0$ the algebra $M_r(A)(t_1, \dots, t_s)$
is defined over a rational extension of $k$.  (Here we are imposing
no restrictions on the degree of $A$.)
\end{enumerate}
\end{thm}

Here $A(t_1, \dots, t_s)$ stands for $A \otimes_F  F(t_1, \dots, t_s)$
and similarly for $M_r(A)(t_1, \dots, t_s)$.
Note that part (b) may be interpreted as saying that for $s \gg 0$
the Brauer class of $A(t_1, \dots, t_s)$ is defined over a rational
extension of $k$.

\subsection{$\bbZ/m \times \bbZ/2$-crossed products}

As usual, we let $(a, b)_m$ denote the symbol
algebra
\begin{equation} \label{E:symbol}
F\{ x, y \}/( x^m = a,  \, y^m = b,  \, xy = \zeta_m yx)\ ,
\end{equation}
where $a,b \in F^*$ and $\zeta_m$ is a (fixed) primitive $m^{\text{th}}$
root
of unity in $F$; cf.~\cite[pp. 194-197]{rowen-rt2}.
In Section~\ref{sect.deg4} we will prove the following:

\begin{thm} \label{thm.ab-cr} Let $A$ be a $\bbZ/m \times \bbZ/2$-crossed
product central simple algebra whose center $F$ contains
a primitive $2m^{\text{th}}$ root of unity $\zeta_{2m}$.
Then $A$ is Brauer equivalent (over $F$) to $(a, b)_m
\otimes_F (c, d)_{2m}$ for some $a, b, c, d \in F^*$.
(In other words,
$M_m(A)$ is isomorphic to $(a, b)_m \otimes_F (c, d)_{2m}$.)
\end{thm}

Note that Theorem~\ref{thm.ab-cr} applies to any
division algebra of degree 4, since, by a theorem of Albert,
any such algebra is a $\bbZ/2 \times \bbZ/2$-crossed product.
In this setting our argument yields, in particular,
an elementary proof of~\cite[Theorem 2, p.~288]{snider}.
We also remark that
Theorem~\ref{thm.ab-cr} may be viewed as an explicit form
of the Merkurjev-Suslin theorem
for $\bbZ/m \times \bbZ/2$-crossed products.

Letting $F_0 = k(\zeta_{2m}, a, b, c, d)$, we note that
\[ (a, b)_m \otimes_F (c, d)_{2m} =
\Bigl( (a, b)_m \otimes_{F_0} (c, d)_{2m}  \Bigr) \otimes_{F_0} F \, .\]
Thus Theorem~\ref{thm.ab-cr} shows that $\Mat_m(A)$ is defined over
$F_0$.  Since $\trdeg_k(F_0) \leq 4$, we obtain the following:

\begin{cor} \label{cor.deg4}
Let $A$ be a  $\bbZ/m \times \bbZ/2$-crossed product central
simple algebra
whose center contains a primitive $2m^{\text{th}}$ root of unity.
Then $\tau(M_m(A)) \le 4$. In particular, $\tau(M_2(A)) \le 4$
for every central simple algebra $A$ of degree 4
whose center contains a primitive $4^{\text{th}}$ root of unity.
\qed
\end{cor}

Note that the last assertion complements, in a somewhat surprising way,
the above-mentioned result of Rost~\eqref{e.rost}. Indeed, suppose
the base field $k$ contains a primitive $4^{\text{th}}$ root of unity.
Then for $A = \UD(4)$, Rost's theorem says that $\tau(A) = 5$,
where as Corollary~\ref{cor.deg4}
says that $\tau(\Mat_2(A)) \leq 4$.

\subsection{Fields of definition of quadratic forms}
\label{subsect1.quadr}

A quadratic form $q \colon V \to F$ on an $F$-vector space
$V = F^n$ is said to be defined over
a subfield $F_0$ of $F$ if $q = q_{F_0} \otimes F$, where $q_{F_0}$
is a quadratic form on $V_0 = F_0^n$.

In the last section we discuss fields of definition
of quadratic forms. Of particular interest to us will be trace forms
of central simple algebras of degree 4, recently studied by Rost,
Serre and Tignol~\cite{rst}. (Recall that the trace form of a central simple
algebra $A$ is the quadratic form $x \mapsto \tr(x^2)$ defined over
the center of $A$.) We will use a theorem of Serre~\cite{serre}
(see our Proposition~\ref{prop.serre}) to prove the following:

\begin{thm} \label{thm.trace-form}
Let $A$ be a central simple algebra of degree $4$ whose
center $F$ contains a primitive $4^{\text{th}}$ root of unity.
Then the trace form of $A$ is defined over a subfield
$F_0 \subset F$ such that $\trdeg_k(F_0) \leq 4$.
\end{thm}

Note that Theorem~\ref{thm.trace-form} may also be viewed
as complementing~\eqref{e.rost}.

\begin{ack} The authors would like to thank the referee for a number
of helpful and constructive comments and for catching several
mistakes in an earlier version of this paper.
\end{ack}

\section{Preliminaries} \label{prelim}

\subsection{$\cG$-lattices}
Throughout, $\cG$ will denote a finite group and $\cH$ will be a subgroup
of $\cG$.
Recall that a \emph{$\cG$-module} is a (left) module over the integral
group ring $\bbZ[\cG]$. As usual, $\Ext_{\cG}$ stands for
$\Ext_{\bbZ[\cG]}$.
A \emph{$\cG$-lattice} is a $\cG$-module
that is free of finite rank over $\bbZ$.
Further, a $\cG$-lattice $M$ is called
\begin{itemize}
\item a \emph{permutation lattice} if $M$ has a
 $\bbZ$-basis that is permuted by $\cG$;
\item \emph{permutation projective} (or \emph{invertible}) if $M$
 is a direct summand of some permutation $\cG$-lattice.
\end{itemize}
A $\cG$-module $M$ is called \emph{faithful}
if no $1\neq g\in\cG$ acts as the identity on $M$.

The $\GH$-augmentation kernel
is defined as the kernel of
the natural augmentation map
$$
\bbZ[\GH]=\bbZ[\cG]\otimes_{\bbZ[\cH]}\bbZ \longrightarrow \bbZ\ ,
\quad \g = g\otimes 1\mapsto 1\qquad (g\in\cG)\ .
$$
Thus there is a short exact sequence
of $\cG$-lattices
\begin{equation}\label{E:augmentation}
0 \to \oGH \to \bbZ[\GH] \to \bbZ \to 0 \ .
\end{equation}
$\omega(\cG/\{1\})$ will be written as $\omega\cG$; this is the
ordinary augmentation
ideal of the group ring $\bbZ[\cG]$; cf.~\cite[Chap.~3]{passman}.

\subsection{Extension sequences}

Exact sequences of $\cG$-lattices of the form
\begin{equation} \label{e.0-1}
\begin{array}{l}
0 \to M \to P \to \oGH \to 0  \, , \\
\text{\rm with $P$ permutation and $M$ faithful}
\end{array}
\end{equation}
will play an important role
in the sequel. In this subsection we introduce two such
sequences,~\eqref{E:freepres} and~\eqref{e.seq2}.

Let $d_{\cG}(\oGH)$ denote the minimum number of
generators of $\oGH$ as a $\bbZ[\cG]$-module.
Then for any $r \geq d_{\cG}(\oGH)$
there exists an exact sequence
\begin{equation} \label{E:freepres}
0 \to M  \to \bbZ[\cG]^r \stackrel{f}{\to} \oGH \to 0\ .
\end{equation}
of $\cG$-lattices.

\begin{lem} \label{L:freepres}
$M$ is a faithful $\cG$-lattice if and only if $r\ge 2$ or $\cH\neq\{1\}$.
\end{lem}

\begin{proof}
It is enough to show that $M\otimes\bbQ$ is $\cG$-faithful;
thus we may work over the semisimple algebra $\bbQ[\cG]$.
Since $f\otimes\bbQ$ splits, we have a $\bbQ[\cG]$-isomorphism
$(\oGH\otimes\bbQ) \oplus (M\otimes\bbQ) \simeq \bbQ[\cG]^r$.
Similarly, the canonical exact sequence
$\bbZ[\cG]\omega\cH\mono\bbZ[\cG]\onto \bbZ[\GH]$ gives
$(\oGH\otimes\bbQ) \oplus \bbQ \oplus \bbQ[\cG]\omega\cH \simeq \bbQ[\cG]$.
Therefore,
$$
M\otimes\bbQ \simeq \bbQ[\cG]^{r-1} \oplus \bbQ \oplus \bbQ[\cG]\omega\cH\ .
$$
If $r\ge 2$ then $\bbQ[\cG]^{r-1}$ is $\cG$-faithful, and
if $\cH \neq \{ 1\}$ then $\omega\cH\otimes\bbQ$ is $\cH$-faithful and so
$\bbQ[\cG]\omega\cH \simeq (\omega\cH\otimes\bbQ)\ind{\cH}{\cG}$ is
$\cG$-faithful.
In either case, $M\otimes\bbQ$ is faithful, as desired. On the other hand,
$r=1$ and $\cH=\{1\}$ leads to $M\otimes\bbQ \simeq \bbQ$ which
is not faithful.
\end{proof}

\begin{lem} \label{permpres}
There is an exact sequence
\begin{equation} \label{e.seq2}
0\to \oGH^{\otimes 2} \stackrel{m}{\to} P \to \oGH\to 0\ ,
\end{equation}
where $P$ is the (permutation) sublattice
$P=\bigoplus_{\x\neq\y\in\GH}\bbZ(\x\otimes{\y})$
of $\bbZ[\GH]^{\otimes 2}$.  The lattice $\oGH^{\otimes 2}$ is faithful if
and only if $\cH$ contains no normal subgroup $\neq\{1\}$ of $\cG$
and $[\cG:\cH]\ge 3$.
\end{lem}

\begin{proof}
Tensoring sequence \eqref{E:augmentation} with $\oGH$ and putting
$P'=\oGH\otimes\bbZ[\GH]$, we obtain an
exact sequence
\[ 0 \to \oGH^{\otimes 2} \to P' \to \oGH \to 0 \ ,
\]
where $\otimes = \otimes_{\bbZ}$.
The elements $(\x-\y)\otimes\y$ with $\x\neq\y\in\GH$ form a $\bbZ$-basis of
$P'$,
and the map
\[ m \colon (\x-\y)\otimes\y\mapsto \x\otimes\y \]
is a $\cG$-isomorphism $P'\simeq P$.

For the faithfulness assertion,
note that $\cN=\bigcap_{g\in\cG}\cH^g$ acts trivially
on $\bbZ[\GH]$ and hence on $\oGH^{\otimes 2}$; so $\cN=\{1\}$ is
surely required
for faithfulness. Also, if $[\cG:\cH]\le 2$ then $\cG$ acts trivially
on $\oGH^{\otimes 2}$. Conversely, if $\cN=\{1\}$ and $[\cG:\cH]\ge 3$ then
it is easy
to verify that $\oGH^{\otimes 2}$ is indeed faithful.
\end{proof}

\subsection{Twisted multiplicative $\cG$-fields} \label{no-name}

Recall that a \emph{$\cG$-field} is a field $F$ on which the finite group
$\cG$ acts by automorphisms, written $f \mapsto g(f)$.
Morphisms of $\cG$-fields are $\cG$-equivariant field
homomorphisms.
The $\cG$-field $F$ is called \emph{faithful} if every $1 \neq g \in \cG$
acts non-trivially on $F$.
If $K \subseteq F$ is a field extension and $V \subseteq F$ a subset of $F$
(not necessarily algebraically independent over $K$)
then we let $K(V)$ denote the subfield of $F$ that is generated
by $K$ and $V$.

\begin{lem} \label{L:no-name} {\rm (cf.~\cite[Appendix 3]{shafarevich})}
Let $K \subseteq F$ be an extension of $\cG$-fields with $K$ faithful.
Assume that $F = K(V)$ for some $\cG$-stable $K$-subspace $V \subseteq F$.
Then
\begin{enumerate}
\item $V = KV^{\cG}$, where $V^{\cG}$ denotes the $\cG$-invariants in $V$,
\item $F = K(V^{\cG})$, and
\item $F^{\cG} = K^{\cG}(V^{\cG})$.
\end{enumerate}
\end{lem}

\begin{proof}
(a) Let $S = K\!\#\cG$ denote the skew group ring for the given $\cG$-action
on
$K$.  The $\cG$-action on $F$ and multiplication with $K$ make $F$ a (left)
$S$-module, and $V$ is a submodule. Moreover, since $K$ is a faithful
$\cG$-field, $S$ is a simple ring; see, e.g., \cite[p.~473]{jacobson}.
In particular, the element $t = \sum_{g \in \cG} g \in S$ generates $S$ as
a 2-sided ideal.
Thus, $S = StS = KtK$ and consequently, $V = KtKV = KV^{\cG}$.

\smallskip
(b) is an immediate consequence of (a).

\smallskip
(c) Let $E = K^{\cG}(V^{\cG})$. We want to show that $E = F^{\cG}$.
Clearly $E \subseteq F^{\cG}$. To prove equality,
note that $KE$ is a subring of $F$ containing $K$ and $V^{\cG}$,
and that $\dim_E KE \le \dim_{K^{\cG}} K = |\cG|$. Thus,
$KE$ is a field, and hence (b) implies that $KE = F$.
Therefore, $\dim_E F \le |\cG| = \dim_{F^{\cG}} F$.
Since $E \subseteq F^{\cG}$, this is only possible
if $E = F^{\cG}$.
\end{proof}

We recall a well-known construction of $\cG$-fields;
cf.~\cite{saltman:twisted}.
Given a $\cG$-field $E$, a $\cG$-lattice $M$, and an extension class
$\gamma\in\Ext_{\cG}(M,E^*)$, the \emph{twisted multiplicative $\cG$-field}
$E_{\gamma}(M)$ is constructed as follows. Form the group algebra $E[M]$ of
$M$ over $E$; this is a commutative integral domain with group of units
$\U(E[M])=E^*\times M$.
We shall use multiplicative notation for $M$ in this setting.
Let $E(M)$ denote the field of fractions of $E[M]$.
Choose an extension of $\cG$-modules
\begin{equation} \label{E:no-name}
1\to E^*\to V\to M\to 1
\end{equation}
representing $\gamma$. So, as abelian groups, $V\simeq \U(E[M])$.
Using this identification,
we obtain a $\cG$-action on $\U(E[M])$ inducing the given action on $E^*$.
The
action of $\cG$ on $\U(E[M])$ extends uniquely to $E[M]$, and to $E(M)$; we
will use $E_{\gamma}[M]$ and $E_{\gamma}(M)$ to denote $E[M]$ and $E(M)$
with
the $\cG$-actions thus obtained. For $\gamma=1$, we will simply write
$E[M]$ and $E(M)$.
We remark that the choice of the sequence~\eqref{E:no-name} representing
a given $\gamma\in\Ext_{\cG}(M,E^*)$ is insubstantial:
a different choice leads to $\cG$-isomorphic results.

For future reference, we record the following application
of Lemma~\ref{L:no-name} essentially
due to Masuda~\cite{masuda}; see also \cite[Proposition 1.6]{lenstra},
\cite[Lemma 12.8]{saltman}.

\begin{prop} \label{P:no-name}
Let $E$ be a faithful $\cG$-field and let $P$ be a permutation
$\cG$-lattice.
Then any twisted multiplicative $\cG$-field $E_{\gamma}(P)$ can be written
as
$$
E_{\gamma}(P)= E(t_1,\ldots,t_n)
$$
with $\cG$-invariant transcendental (over $E$) elements $t_i$. In
particular,
$E_{\gamma}(P)^{\cG}=E^{\cG}(t_1,\ldots,t_n)$ is rational over $E^{\cG}$.
\end{prop}

\begin{proof}
We have an extension of $\cG$-modules
$1 \to E^*\to \U(E_\gamma[P]) \to P \to 1$
representing $\gamma$, as in \eqref{E:no-name}.
Fix a $\bbZ$-basis, $X$, of $P$ that is permuted by the action of $\cG$.
For each $x \in X$, choose a preimage
$x' \in \U(E_\gamma[P])$.
Then $\{x'\}_{x \in X}$ is
a collection of transcendental generators of $E_\gamma(P)$ over $E$, and
$\cG$ acts via
$g(x') = g(x)'y$ for some $y = y(g,x) \in E^*$.
Letting $V$ denote the $E$-subspace of
$E_\gamma(P)$ that is generated by
$\{x'\}_{x \in X}$, we conclude from Lemma~\ref{L:no-name} that
$V$ has a basis consisting of $\cG$-invariant elements, say
$t_1,\ldots,t_n$,
and $E_{\gamma}(P)= E(t_1,\ldots,t_n)$,
$E_{\gamma}(P)^{\cG}=E^{\cG}(t_1,\ldots,t_n)$.
The $t_i$ are transcendental over $E$, since $\trdeg_E E_\gamma(P) =
\rank(P)
= n$.
\end{proof}


\subsection{Rational specialization}
\label{sect2.3}

Let $A/F$ and $B/K$ be central simple algebras.
We will call $B/K$ a \emph{rational specialization} of $A/F$ if there exists
a field $F'$ containing both $F$ and $K$ such that $F'/K$ is rational
and \[ B \otimes_K F' \simeq A\otimes_F F' \, . \]
In other words, $B$ is a rational specialization of $A$ if
$\deg A = \deg B$ and $A$ embeds in some $B(t_1, \dots, t_n)$,
where $t_1, t_2, \dots$ are independent variables over $F$;
cf.~\cite[p.~73]{saltman},

\smallskip

\emph{For the rest of this paper we fix an (arbitrary) base field $k$.
All other fields are understood to contain a copy $k$ and all
maps (i.e., inclusions) between fields are understood to restrict
to the identity map on $k$.}

\begin{defn} \label{def.rat-spec}
Let $\Lambda$ be a class of central simple algebras. We shall say
that an algebra $A \in \Lambda$ has the \emph{rational specialization
property} in
the class $\Lambda$ if every $B \in \Lambda$ is a rational specialization
of $A$. If $\Lambda$ is the class of all central simple algebras of degree
$n = \deg(A)$ then we will omit the reference to $\Lambda$
and will simply say that $A$ has the rational specialization property.
\end{defn}

\begin{example} \label{ex.rat-spec}
By \cite[Lemma 3.1]{rv2}, $\UD(n)$ has the rational specialization
property.  This is also implicit in~\cite{saltman-note}.
We remark that any central simple algebra $A/F$ with the rational
specialization
property is a division algebra. To see this, specialize $A$ to $\UD(n)$,
where $n = \deg(A)$.
\end{example}

Recall the definition of $\tau(A)$ given at the beginning of this paper.

\begin{lem}\label{rattau}
Let $A/F$ and $B/K$ be a central simple algebras.
\begin{enumerate}
\item If $A'\simeq A \otimes_{F} F'$ for some rational
field extension $F'/F$ then $\tau(A)=\tau(A')$.

\item \textnormal{(cf.~\cite[Lemma 11.1]{saltman})}
If $A$ is a rational specialization of $B$ then
$\tau(A)\le \tau(B)$.
\end{enumerate}
\end{lem}

\begin{proof}
(a) The inequality $\tau(A') \le \tau(A)$ is immediate from the definition
of $\tau$. To prove the opposite inequality,  suppose
$A'\simeq A_0 \otimes_{F_0} F'$, where $A_0$ is a central simple algebra
over
an intermediate field $k \subset F_0 \subset F'$,
and $A_0$ is chosen so that $\trdeg_k(F_0) = \tau(A')$.
In particular, $\trdeg_k(F_0) \leq \trdeg_k(F)$.
Then by~\cite[Proposition 3.2]{rv2}, $A_0$ embeds in $A$, i.e.,
$A \simeq A_0 \otimes_{F_0} F$ for some embedding $F_0 \hookrightarrow F$.
Consequently, $\tau(A) \le \trdeg_k(K) = \tau(A')$, as desired.

\smallskip
(b) We may assume $B \otimes_{K} F' = A'$, as in (a).
Clearly $\tau(B) \ge \tau(A')$, and part (a) tells us that
$\tau(A') = \tau(A)$.
\end{proof}

\begin{remark} \label{rem-ud(n)}
Combining Example~\ref{ex.rat-spec} with
Lemma~\ref{rattau}(b), we see that
$\tau(A) \leq d(n)=\tau(\UD(n))$ holds for every central simple algebra $A$
of degree $n$; cf.~\cite[Lemma 9.2]{reichstein1}.
\end{remark}

\section{$\GH$-crossed products}
\label{crossed}

We shall call a central simple algebra $A/F$ an \emph{$(E,\GH)$-crossed
product}
if $A$ has a maximal subfield $L$ whose Galois closure $E$ over $F$
has the property that $\Gal(E/F) = \cG$ and $\Gal(E/L) = \cH$.
(We adopt the convention that a \emph{maximal subfield} of $A$ is a subfield
$L$ that is maximal as a commutative subring; so $[L : F]$ is equal to the
degree of $A$.)
We will say that $A$ is a $\GH$-crossed product if it is an
$(E,\GH)$-crossed product for some faithful $\cG$-field $E$.
If $\cH = \{1 \}$ then a $\GH$-crossed product is just a
\emph{$\cG$-crossed product} in the usual sense
(see, e.g., \cite[Definition 3.1.23]{rowen-pi}).

\begin{example} \label{ex.procesi}
Consider the universal division algebra $\UD(n)$
generated by two generic matrices, $X$ and $Y$, over $k$. Denote the
center of this algebra by $Z(n)$. Setting $L = Z(n)(X)$,
we see that $\UD(n)$ is an $\Sym_n/\Sym_{n-1}$-crossed
product \cite[Theorem 1.9]{procesi}; see also Section~\ref{sect.ud} below.
\end{example}

Since the degree of a $\GH$-crossed product is equal to $[\cG : \cH]$,
we see that isomorphism classes of $(E,\GH)$-crossed products are
in 1-1 correspondence with
the relative Brauer group $\BB(L/F)$, which, in turn, is naturally
identified with the kernel of the restriction homomorphism
$H^2(\cG, E^*) \to H^2(\cH, E^*)$; cf. \cite[14.7]{pierce}.

A $\cG$-module $M$ is called
\emph{$H^1$-trivial} if $H^1(\cH,M)=0$ holds for every subgroup
 $\cH\le\cG$. Equivalently, $M$ is $H^1$-trivial
if $\Ext_{\cG}(P,M)=0$ for all permutation projective $\cG$-lattices $P$;
see, e.g., \cite[Lemma 12.3]{saltman}.

\begin{lem} \label{lem.tue-1}
Given an exact sequence
\[ 0 \to M \to P \to \oGH \to 0  \, , \]
of $\cG$-lattices, with $P$-permutation,
let $N$ be an $H^1$-trivial $\cG$-module.
Denote the kernel of the restriction
homomorphism $H^2(\cG, N) \to H^2(\cH, N)$ by
$K(\GH, N)$.  Then there is a natural isomorphism
$$
\phi_N \colon
\Hom_{\cG}(M, N)/\im(\Hom_{\cG}(P, N)) \stackrel{\simeq}{\longrightarrow}
K(\GH, N)\ .
$$
\end{lem}

Here ``natural" means
that for every homomorphism $N \to N'$
of $H^1$-trivial $\cG$-modules, the following diagram commutes
\begin{equation} \label{e.trivial}
\xymatrix{
\Hom_{\cG}(M, N')/\im(\Hom_{\cG}(P, N'))  \ar[r]^-{\phi_{N'}} & K(\GH, N')
\\
\Hom_{\cG}(M, N)/\im(\Hom_{\cG}(P, N)) \ar[u]  \ar[r]^-{\phi_N} & K(\GH, N)
\ar[u]
}
\end{equation}
Note that, other than in sequence \eqref{e.0-1},
the $\cG$-lattice $M$ need not be faithful.

\begin{proof} The lemma is a variant of~\cite[Theorem 12.10]{saltman},
where the same assertion is made for the sequence~\eqref{e.seq2}.
The proof of~\cite[Theorem 12.10]{saltman}
goes through unchanged in our setting.
\end{proof}

In subsequent applications we will always take $N = E^*$, where $E$ is
a faithful $\cG$-field. Note that $E^*$ is $H^1$-trivial
by Hilbert's Theorem 90.
As we remarked above, $K(\GH, E^*)$ is naturally identified with
$\Br(L/F)$, where $L = E^{\cH}$, and elements of $\Br(L/F)$
are in 1-1 correspondence with $(E, \GH)$-crossed products.
We shall denote the $(E, \GH)$-crossed product associated to
a $\cG$-homomorphism $f \colon M \to E^*$ by $\Alg(f)$.

\begin{lem} \label{lem.tue-2}
Consider a sequence of $\cG$-lattices of the form~\eqref{e.0-1}.
Let $E$ be $\cG$-field and $f \colon M \to E^*$ be a homomorphism
of $\cG$-modules.
If $k(f(M))$ is contained in a faithful $\cG$-subfield $E_0$ of $E$ then
$\Alg(f)$ is defined over $E_0^{\cG}$.
\end{lem}

\begin{proof} Since $f$ is
the composition of $f_0 \colon M \to E_0^*$ with the inclusion
$E_0^* \into E^*$, Lemma~\ref{lem.tue-1} tells us that
$A = \Alg(f_0) \otimes_{E_0^{\cG}} E^{\cG}$.
\end{proof}

\begin{remark} \label{rem.bfs} In the special case where
the sequence \[ 0 \to M \to P \to \oGH \to 0 \] is given by~\eqref{e.seq2},
$M = \oGH^{\otimes 2}$ has a particularly convenient
set of generators
\[ y_{ijh} = (\overline{g_i} - \overline{g_j}) \otimes
(\overline{g_j} - \overline{g_h}) \, , \]
where $\cG / \cH = \{ \overline{g_1}, \dots,
\overline{g_n} \}$ is the set of left cosets of $\cH$ in $\cG$ and
$i$, $j$, $h$ range from $1$ to $n = [\cG: \cH]$;
cf.~\cite[Lemma 1.2]{rowensaltman}.
If $f \colon \oGH^{\otimes 2} \to E^*$ is a $\cG$-module
homomorphism then the elements $c_{ijh} = f(y_{ijh})$
form a {\em reduced Brauer factor set} for $\Alg(f)$ in the sense of
\cite[p. 449]{rowensaltman}.
Conversely, for any reduced Brauer factor set
$(c_{ijh})$ in $E^*$,
there exists a homomorphism $f \colon \oGH^{\otimes 2}
\to E^*$ such that $f(y_{ijh}) = c_{ijh}$;
see~\cite[Corollary 1.3]{rowensaltman}. Thus
Lemma~\ref{lem.tue-2} takes the following form:

\begin{quote}
\emph{Let $A$ be an $(E, \GH)$-crossed product defined by a reduced Brauer
factor set $(c_{ijh})$. Suppose $(c_{ijh})$ is contained in a faithful
$\cG$-subfield $E_0$ of $E$. Then
$A$ is defined over $E_0^{\cG}$.}
\qed
\end{quote}
\end{remark}

This following theorem is a variant of~\cite[Theorem 12.11]{saltman}.

\begin{thm} \label{thm.tue-3}
Given the sequence~\eqref{e.0-1}, let
$\mu: M \hookrightarrow k(M)^*$ be the natural inclusion.
Then $D= \Alg(\mu)$ has the rational specialization property in the class
of $\GH$-crossed products containing a copy of $k$ in their center.
In particular, $\tau(A) \leq \rank(M)$ holds for
any $\GH$-crossed product $A/F$ with $k \subset F$.
\end{thm}

\begin{proof} Write $A = \Alg(f)$ for some
$\cG$-homomorphism $f \colon M \to E^*$, where $E$ is a faithful
$\cG$-field with $E^{\cG}=F$;
see the remarks following
Lemma~\ref{lem.tue-1}. Furthermore, let
$E(P)$ denote the fraction field of the group algebra $E[P]$, with the
$\cG$-action induced from the $\cG$-actions on $E$ and $P$.
By Proposition~\ref{P:no-name}, there exists an $E$-isomorphism
$j \colon E(P) \simeq E(\bt)$ of $\cG$-fields,
where $\bt = (t_1, \dots, t_r)$, $r = \rank(P)$, are
indeterminates on which $\cG$ acts trivially.
Therefore, $E(P)^{\cG} \simeq E^{\cG}(\bt) = F(\bt)$ is a rational
extension of $F$.
Let $f_{\bt} \colon M \to E(\bt)^*$ denote the composition of
$f$ with the natural inclusion $E^* \into E(\bt)^*$. Then
$\Alg(f_{\bt}) = \Alg(f)\otimes_F F(\bt) = A\otimes_F F(\bt)$.
By Lemma~\ref{lem.tue-1}, $\Alg(f_{\bt}) \simeq \Alg(f_{\bt}+g|_M)$ for
any $g \in \Hom_{\cG}(P, E(\bt)^*)$.  Let $g$ be the composite
$g \colon P \into E(P)^* \stackrel{\sim}{\longrightarrow} E(\bt)^*$
and let $\varphi$ be the $\cG$-module map $\varphi \colon M \to E(\bt)^*$,
$\varphi(m)=f_{\bt}(m)g(m)$. We claim that $\varphi$
lifts to an embedding of $\cG$-fields $k(M) \into E(\bt)$.
Indeed, modulo $E^*$, $\varphi(m) \equiv g(m) \in P \subseteq E(\bt)^*$.
Hence,
$\{\varphi(m)\}_{m \in M}$ is an $E$-linearly independent subset of
$E(\bt)$,
and so the map $k[\varphi] \colon k[M] \to E(\bt)$,
$\sum_m k_mm \mapsto \sum_m k_m\varphi(m)$ is a
$\cG$-equivariant embedding of the group ring $k[M]$ into $E(\bt)$.
This embedding lifts to an embedding of
$\cG$-fields $\phi \colon k(M) = Q(k[M]) \into E(\bt)$, as we have claimed.
So $\phi\circ\mu = \varphi$, and hence
$D\otimes_{k(M)^{\cG}}F(\bt) = \Alg(\phi\circ\mu) = \Alg(\varphi) \simeq
\Alg(f_{\bt}) = A\otimes_F F(\bt)$. This proves that $A$ is a rational
specialization of $D$.

Lemmas~\ref{rattau}(b) and \ref{lem.tue-2} now imply that
$\tau(A) \leq \tau(D) \leq \trdeg_k k(M)^{\cG} = \rank(M)$.
This completes the proof of the theorem.
\end{proof}

\begin{remark} \label{rem.exponent}
Continuing with the notation used in the above theorem,
the rational specialization property of $D=\Alg(\mu)$
implies that \emph{$D$ is a division algebra of
exponent $[\cG:\cH]$}. Indeed, by~\cite[Appendix]{fss} there
exists a $\GH$-crossed product division algebra of exponent
$[\cG:\cH]$, and the above assertion can be proved by specializing
$D$ to this algebra.  Alternatively,
the fact that $D$ is a division algebra of exponent
$[\cG:\cH]$ can be checked directly by showing that the image
of $\mu$ in $H^2(\cG,k(M)^*)$ (see Lemma~\ref{lem.tue-1}) has order
$[\cG:\cH]$.
\end{remark}

\begin{remark} \label{rem.freepres}
The above construction applies in particular to sequences of
the form~\eqref{E:freepres}. The following special
type of sequence~\eqref{E:freepres}
has been particularly well-explored.
Write $\cG=\langle\cH,g_1,\dots,g_r\rangle$
for suitable $g_i\in\cG$;
the minimal such $r$ is usually denoted by $d(\GH)$.
Then we can define an epimorphism of
$\cG$-lattices
$f \colon \bbZ[\cG]^r \onto \oGH$, $f(\alpha_1,\dots,\alpha_r)=
 \sum_{i=1}^r\alpha_i\overline{(g_i-1)}$, where
$\overline{\phantom{iii}}\colon \bbZ[\cG]\onto \bbZ[\GH]$ is the canonical
map;
see~\cite[Lemma 3.1.1]{passman}.
The kernel $R(\GH)=\Ker f$ is called a \emph{relative relation
module}; it has the following group theoretical description. Let
$\cF_r$ denote the free group on $r$ generators and consider the
presentation
\begin{equation} \label{E:relpres}
1 \to \cR \to \cF_r\ast\cH \to \cG \to 1 \,
\end{equation}
where $\cF_r\ast\cH \to \cG$ is the identity on $\cH$ and sends the $r$
generators of
$\cF_r$ to the elements $g_1,\dots,g_r$.
Then $R(\GH) \simeq \cR^{\text{ab}}=\cR/[\cR,\cR]$,
with $\cG$ acting by conjugation; see \cite{kimmerle} and \cite{gruenberg}
(for $\cH=\{1\}$).
Thus, we have the following version of sequence~\eqref{E:freepres}
with $M=\cR^{\text{ab}}$ :
\begin{equation} \label{E:relationmod}
0 \to \cR^{\text{ab}} \to \bbZ[\cG]^r \to \oGH \to 0 \ .
\end{equation}
When $\cH=\{1\}$ and $r\ge 2$, the division algebra $D$ constructed via
\eqref{E:relationmod} in Theorem~\ref{thm.tue-3} is identical with
the generic $\cG$-crossed product of Snider~\cite{snider};
see also Rosset~\cite{rosset}.
Explicitly:
\begin{quote}
\emph{Given a free presentation $1\to\cR\to\cF_r\to\cG\to 1$ of $\cG$
with $r\ge 2$,
let $M=\cR^{\text{\rm ab}}\le \overline{\cF_r} = \cF_r/[\cR,\cR]$ and
let $a\in H^2(\cG,k(M)^*)$ be the image of the extension class
$[1\to M \to\overline{\cF_r} \to\cG\to 1] \in H^2(\cG,M)$ under
the natural inclusion $\mu: M \hookrightarrow k(M)^*$.
Then $D=\Alg(\mu)$ is the $\cG$-crossed product $(k(M),\cG,a)$
or, equivalently,
the localization of the group algebra $k[\overline{\cF_r}]$
at the nonzero elements of $k[M]$.}
\end{quote}
\end{remark}

\begin{cor} \label{cor2-generic}
Let $A$ be a $\GH$-crossed product
and let $d_{\cG}(\oGH)$ be the minimal number of generators
of $\oGH$ as a $\cG$-module.
Then \[ \tau(A) \le r|\cG|-[\cG:\cH]+1 \, , \] where
\[ r =
\begin{cases}
d_{\cG}(\oGH) & \text{if $\cH\neq\{1\}$} \\
\max\{2,d_{\cG}(\oGH)\} & \text{if $\cH = \{1\}$}
\end{cases} \]
\end{cor}

\begin{proof} Applying Theorem~\ref{thm.tue-3} to the exact
sequence~\eqref{E:freepres}, we obtain
\[ \tau(A) \leq \rank(M) = \rank(\bbZ[G]^r) - \rank(\oGH) =
 r|\cG|-[\cG:\cH]+1 \, , \]
as claimed.
Note that for $r$ as above, Lemma~\ref{L:freepres} tells us that $M$
is faithful, so that Theorem~\ref{thm.tue-3} is, indeed, applicable.
\end{proof}

\begin{remark}
As we pointed out in Remark~\ref{rem.freepres},
$d_{\cG}(\oGH) \le d(\GH)$.
The difference $\pr(\GH)=d(\GH) - d_{\cG}(\oGH) \ge 0$
can be arbitrarily large, even if $\cH=\{1\}$.
In this case $d(\cG)=d(\cG/\{1\})$ is the minimal
number of generators of $\cG$, and
$\pr(\cG)=d(\cG)-d_{\cG}(\omega\cG)$ is usually called
the \emph{presentation rank} or \emph{generation gap} of $\cG$.
All solvable groups $\cG$ have presentation rank $\pr(\cG)=0$;
see~\cite[Lectures 6 and 7]{gruenberg}.
Moreover, if the derived subgroup $[\cG,\cG]$ is nilpotent then
$\pr(\GH)=0$ holds for every subgroup $\cH$ of $\cG$; see~\cite{kimmerle}.
\end{remark}

\begin{cor} \label{cor2a}
\begin{enumerate}
\item Suppose a group $\cG$ of order $n$ can be generated
by $r\ge 2$ elements. Then
$\tau(A) \leq (r-1)n + 1$ for any $\cG$-crossed product central
simple algebra $A$.
\item $\tau(A) \leq (\lfloor \log_2(n) \rfloor - 1)n + 1$ holds for
any crossed product central simple algebra $A$ of degree $n \geq 4$.
\end{enumerate}
\end{cor}
Here, as usual, $\lfloor x \rfloor$ denotes the largest integer $\leq x$.

\begin{proof} (a) is an immediate
consequence Corollary~\ref{cor2-generic}. (b) follows from
(a), because any group of order $n$ can be generated
by $r  \leq \log_2(n)$ elements.
(Indeed, $|\langle \cG_0, g \rangle| \geq 2 |\cG_0|$
for any subgroup $\cG_0$ of $\cG$ and any $g \in \cG \setminus \cG_0$.)
Note also that $\lfloor \log_2(n) \rfloor \geq 2$ for any $n \ge 4$.
\end{proof}

\section{Proof of Theorem~\ref{thm1}}
\label{sect.ud}

For the next two sections we shall assume that $\cG = \Sym_n$ and
$\cH = \Sym_{n-1}$. We will use the following standard notations
for $\Sym_n$-lattices:
\begin{equation} \label{e.s_n-lat}
 \bbZ[\Sym_n/\Sym_{n-1}] = U_n \quad \text{\rm and} \quad
\omega(\Sym_n/\Sym_{n-1}) = A_{n-1} \, .
\end{equation}
The natural generators of $U_n$ will be denoted by $u_1, \dots, u_n$;
the symmetric group $\Sym_n$ permutes them via $\sigma(u_i) =
u_{\sigma(i)}$.
$A_{n-1}$ is the sublattice of $U_n$ generated by $u_i - u_1$ as $i$
ranges from $2$ to $n$.

Recall that the universal division algebra $\UD(n)$
is generated, as a $k$-division algebra, by
a pair of generic $n \times n$-matrices $X$ and $Y$.
We may assume without loss of generality that $X$ is diagonal.
Following \cite{rowen-brauer} we will denote
the diagonal entries of $X$ by $\zeta'_{ii}$ and
the entries of $Y$ by $\zeta_{ij}$, where $\zeta'_{ii}$ and
$\zeta_{ij}$ are algebraically independent variables over $k$.
The group $\Sym_n$ permutes these variables as follows:
\[ \text{$\sigma(\zeta'_{ii}) = \zeta'_{\sigma(i) \sigma(i)}$ and
$\sigma(\zeta_{ij}) = \zeta_{\sigma(i) \sigma(j)}$.} \]
We identify the multiplicative group generated by $\zeta_{ii}'$
with the $\Sym_n$-lattice $U_n$ (via $\zeta_{ii}' \leftrightarrow u_i$),
and the multiplicative group generated
by $\zeta_{ij}$ with $U_n \otimes U_n$
(via $\zeta_{ij} \leftrightarrow u_i \otimes u_j$).
Consider the exact sequence
\begin{equation} \label{e.formanek-procesi}
0 \to \Ker(f) \to U_n \oplus U_n^{\otimes 2}
\stackrel{f}{\to} A_{n-1} \to 0
\end{equation}
of $\Sym_n$-lattices, where $f(u_i, u_j \otimes u_h) = u_j - u_h$.
This sequence is the sequence \eqref{e.seq2} of
Lemma~\ref{permpres} for $\cG = \Sym_n$ and $\cH = \Sym_{n-1}$, with two
extra copies of $U_n$ added: the second copy
of $U_n$ is the sublattice of $U_n^{\otimes 2}$ that is spanned by all
elements
$u_i \otimes u_i$. Both copies of $U_n$ belong to $\Ker(f)$; in fact,
$$
\Ker(f) = U_n \oplus U_n \oplus A_{n-1}^{\otimes 2} \ ,
$$
where $A_{n-1}^{\otimes 2}$ is identified with the sublattice of
$U_n^{\otimes 2}$
that is spanned by all elements $(u_i - u_j) \otimes (u_l - u_m)$.

Let $E = k(\Ker(f))$ and $F = E^{\Sym_n}$. By a theorem
of Formanek and Procesi, $F$ is naturally isomorphic
to the center $Z(n)$ of $\UD(n)$; see,
e.g.,~\cite[Theorem 3]{formanek-3x3}.
Note that $E = F(\zeta_{11}', \dots, \zeta_{nn}')$ is generated over $F$
by the eigenvalues of the generic matrix $X$. Consequently,
$\UD(n)$ is an $(E, \Sym_n/\Sym_{n-1})$-product, and
$E^{\Sym_{n-1}}$ is isomorphic to the maximal subfield $Z(n)(X)$ of
$\UD(n)$;
see,~\cite[Section II.1]{procesi}.

Theorem~\ref{thm1} is now a consequence of the following:

\begin{prop} \label{prop.bfs2}
Suppose $n \geq 5$ is odd. Then
\begin{enumerate}
\item $\UD(n)$ is defined over
$F_0 = k(\bigwedge^2 A_{n-1})^{\Sym_n}$,
\item $Z(n) = k(\Ker(f))^{\Sym_n}$ is rational over
$F_0 = k(\bigwedge^2 A_{n-1})^{\Sym_n}$,
\end{enumerate}
\end{prop}

Here, we view $\wed{A_{n-1}}$ as the sublattice of antisymmetric tensors in
$A_{n-1}^{\otimes 2}$, that is, the $\bbZ$-span of all
$a \wedge a' = a \otimes a' - a' \otimes a$ with $a,a' \in A_{n-1}$.

\begin{proof}
We will deduce part (a) from Remark~\ref{rem.bfs}
by constructing a reduced Brauer factor
set contained in $E_0 = k(\wed{A_{n-1}})$. First we note that
the $\Sym_n$-action on $E_0$ is faithful, because $\wed{A_{n-1}}$
is a faithful $\Sym_n$-lattice for every $n \geq 4$. (Indeed,
$\wed{A_{n-1}} \otimes\bbQ$
is the simple $\Sym_n$-representation corresponding to the
partition $(n-2,1^2)$ of $n$; cf.~\cite[Exercise 4.6]{fulton}).

We now proceed with the construction of the desired Brauer factor set.
The computation in~\cite[Section 2]{rowen-brauer}
shows that the elements
\[ c_{ijh} = \zeta_{ij} \zeta_{jh} \zeta_{ih}^{-1} \in E^*
\, . \]
form a Brauer factor set for $\UD(n)$.
By~\cite[Theorem 4]{rowen-brauer}, if $n$ is odd, $\UD(n)$
has a normalized (and, in particular, reduced)
Brauer factor set $(c'_{ijh})$ given by
\[ c'_{ijh} = (c_{ijh}/c_{hji})^{\frac{n+1}{2}} = (\zeta_{ij}
\zeta_{ji}^{-1}
\zeta_{jh} \zeta_{hj}^{-1} \zeta_{hi} \zeta_{ih}^{-1})^{\frac{n+1}{2}}
\, . \]
Now observe that $\zeta_{ij} \zeta_{ji}^{-1}
\zeta_{jh} \zeta_{hj}^{-1} \zeta_{hi} \zeta_{ih}^{-1}$
is precisely the element of $U_n^{\otimes 2}$ we identified
with $(u_i - u_j) \wedge (u_j - u_h)$.  Thus every $c'_{ijh}$ lies in
$\bigwedge^2{A_{n-1}} \subset E_0$, as desired.

\smallskip
(b) The canonical exact sequence
\[ 0 \to \wed{A_{n-1}} \tto A_{n-1}^{\otimes 2} \tto \sym{A_{n-1}} \to 0 \]
of $\Sym_n$-lattices gives rise to an exact sequence
\[ 0\to\wed{A_{n-1}}
\to A_{n-1}^{\otimes 2}\oplus U_n\oplus \bbZ \to Q\to 0 \ ,
\]
where we have put $Q=\sym{A_{n-1}}\oplus U_n\oplus\bbZ$.  The crucial
fact here is that, by \cite[Section 3.5]{LL}, if $n$ is odd,
$Q$ is a permutation lattice.
Applying Proposition~\ref{P:no-name} to the extension of (faithful)
$\Sym_n$-fields
$E_0 = k(\wed{A_{n-1}}) \subseteq k(A_{n-1}^{\otimes 2}\oplus U_n\oplus\bbZ)
\cong (E_0)_\gamma(Q)$, where $\gamma$ is the image of the class of the
above
extension in $\Ext_\cG(Q,E_0^*)$,
we conclude that
$$ \textstyle
k(A_{n-1}^{\otimes 2}\oplus U_n\oplus\bbZ) \simeq
k(\bigwedge^2 \, A_{n-1})(x_1,\ldots,x_m)
$$
as $\Sym_n$-fields, where $m = \frac{n(n+1)}{2} + 1$ and
$\Sym_n$ acts trivially on the $x_i$'s.  Similarly,
putting $L_n=A_{n-1}^{\otimes 2} \oplus U_n^2$, the obvious
sequence
$0\to A_{n-1}^{\otimes 2} \oplus U_n \to L_n\to U_n\to 0$ leads to
$
k(L_n) \simeq k(A_{n-1}^{\otimes 2} \oplus U_n)(t_1,\ldots,t_n)
$
as $\Sym_n$-fields.
Therefore,
\begin{eqnarray*}
\textstyle
k(L_n) &\simeq&
k(A_{n-1}^{\otimes 2} \oplus U_n)(t_1,\ldots,t_n)\\
&=& k(A_{n-1}^{\otimes 2}\oplus U_n\oplus\bbZ)(t_1,\ldots,t_{n-1})\\
&\simeq& k(\wed{A_{n-1}})(x_1,\ldots,x_m,t_1,\ldots,t_{n-1})
\end{eqnarray*}
as $\Sym_n$-fields, which implies that
$$
Z(n) \simeq k(L_n)^{\Sym_n} \simeq
k(\wed{A_{n-1}})^{\Sym_n}(x_1,\ldots,x_m,t_1,\ldots,t_{n-1})\ ;
$$
so $Z(n)$ is rational over $F_0 = k(\wed{A_{n-1}})^{\Sym_n}$.
\end{proof}

\section{Proof of Theorem~\ref{thm.rat}}
\label{sect.rat}

\subsection{Proof of part (a)}

\begin{reduction} \label{red5.1} Suppose an algebra $A_0$ of degree $n$
has the rational specialization property (see Section~\ref{sect2.3}).
If Theorem~\ref{thm.rat}(a) holds for $A_0$ then it holds for any
central simple algebra $A$ of degree $n$.
\end{reduction}

\begin{proof}  Suppose the for some $r \geq 1$, $A_0(t_1, \dots, t_r)$
is defined over a rational extension $F_0$ of $k$.
Let $A/F$ be an arbitrary central simple algebra of degree $n$.
Then by the rational specialization property, $A_0$ embeds in
$A(t_{r+1}, \dots, t_s)$ for some $s \gg 0$; thus $A_0(t_1, \dots, t_r)$
embeds in $A(t_1, \dots, t_s)$. This shows that $A(t_1, \dots, t_s)$
is defined over $F_0$, as desired.
\end{proof}

In particular, in proving Theorem~\ref{thm.rat}(a), we may assume that
$A$ is a division algebra of degree $n$; see Example~\ref{ex.rat-spec}. By
primary decomposition (cf., e.g., \cite[p.~261]{pierce}),
we only need to consider the cases
where $n = 2$, $n = 4$ and $n$ is an odd prime.
This follows from the next reduction:

\begin{reduction} \label{red5.2} If the conclusion of
Theorem~\ref{thm.rat}(a)
holds for central simple algebras $A_1/F$ and $A_2/F$ (for every choice of
the base field $k \subset F$) then it also holds for $A = A_1 \otimes_F
A_2$.
\end{reduction}

\begin{proof} After replacing $A_1$ and $A_2$ by, respectively,
$A_1(t_1, \dots, t_s)$ and
$A_2(t_1, \dots, t_s)$, we may assume that $A_1$ is defined over
a subfield $F_1 \subset F$ such that $k \subset F_1$ and
$F_1$ is rational over $k$. We will now think of $F_1$ (rather than $k$)
as our new base field. After adding more indeterminates, we may
assume that $A_2$ is defined over a subfield $F_2 \subset F$,
where $F_1 \subset F_2$ and $F_2$ is rational over $F_1$.
Now $F_2$ is rational over $k$, and since $A_1$ and $A_2$ are both
defined over $F_2$, so is $A$.
\end{proof}

We are now ready to complete the proof of Theorem~\ref{thm.rat}(a).

First, suppose $n = 2$ or $4$. Since $\UD(n)$ has the rational
specialization property, we may assume $A = \UD(n)$; see
Reduction~\ref{red5.1}. But since the center of $\UD(n)$ is known
to be rational for $n = 2$ (see~\cite[Theorem 2.2]{procesi})
and $n = 4$ (see~\cite{formanek-4x4}), these algebras clearly
satisfy the conclusion of Theorem~\ref{thm.rat}(a). This completes
the proof of the theorem for $n = 2$ and $4$.
We remark that the same argument goes through
for $n = 3$ (because the center of $\UD(3)$ is known to be rational;
see \cite{formanek-3x3}) and for $n = 5, 7$ (because the centers
of $\UD(5)$ and $\UD(7)$ are known to be stably rational; see~\cite{bl}),
but we shall not need it in these cases.

From now on we will assume that $n = p$ is an odd prime.
Then the $\Sym_n$-lattice $A_{n-1}^{\otimes 2}$ is faithful;
see Lemma~\ref{permpres}.
Furthermore, by a theorem of Bessenrodt and
LeBruyn~\cite[Proposition 3]{bl} (see also \cite[Lemma 2.8]{beneish} for a
more
explicit form of this result),  $A_{n-1}^{\otimes 2}$
is permutation projective, i.e., there exists an $\Sym_n$-lattice $L$
such that $P = A_{n-1}^{\otimes 2} \oplus L$ is permutation.
We can assume that $k(P)^{\Sym_n}$ is rational over $k$. Indeed, after
adding a
copy of $U_n$ if necessary, we have $P = U_n \oplus Q$ for some permutation
lattice $Q$, and so $k(P) \simeq k(U_n)(Q)$. Proposition~\ref{P:no-name}
implies that $k(P)^{\Sym_n}$ is rational over $k(U_n)^{\Sym_n}$,
which in turn is rational over $k$.

Let
\[ i \colon A_{n-1}^{\otimes 2} \hookrightarrow
k(A_{n-1}^{\otimes 2})^*  \]
and
\[ j \colon A_{n-1}^{\otimes 2} \hookrightarrow
k(A_{n-1}^{\otimes 2} \oplus P)^*
\, \]
be the natural embeddings of $\Sym_n$-modules. (Here, $j$ identifies
$A_{n-1}^{\otimes 2}$ with the first summand of
$A_{n-1}^{\otimes 2} \oplus P$.)
Recall that by Lemma~\ref{lem.tue-1}
these embeddings, in combination with the exact sequence~\eqref{e.seq2}
(for $\cG = \Sym_n$ and $\cH = \Sym_{n-1}$; see ~\eqref{e.s_n-lat}),
give rise to central simple algebras $\Alg(i)$ and $\Alg(j)$.

By Theorem~\ref{thm.tue-3}, $\Alg(i)$ has the rational specialization
property
in the class of $\Sym_n/\Sym_{n-1}$-crossed
products. Thus,
the universal division algebra $\UD(n)$, being an
$\Sym_n/\Sym_{n-1}$-crossed
product (see Example~\ref{ex.procesi}), is a rational
specialization of $\Alg(i)$.
Since $\UD(n)$ has the rational specialization property
in the class of all central simple algebras of degree $n$ (see
Example~\ref{ex.rat-spec}), so does $\Alg(i)$.

We claim that $\Alg(j)$ also has the rational specialization property
in the class of central simple algebras of degree $n$.
Indeed, by Lemma~\ref{lem.tue-1},
\[ \Alg(j) = \Alg(i) \otimes_{F^{\Sym_n}} E^{\Sym_n} \, , \]
where $F = k(A_{n-1}^{\otimes 2})$ and $E =
k(A_{n-1}^{\otimes 2} \oplus P)$.
Now Proposition~\ref{P:no-name} tells us that $E^{\Sym_n}$ is a rational
extension of $F^{\Sym_n}$, and the claim follows.

By Reduction~\ref{red5.1} it now suffices to prove that
$\Alg(j)$ is defined over a purely transcendental extension of $k$.
Put $E_0 = k(A_{n-1}^{\otimes 2} \oplus (0) \oplus L) \subseteq E$.
Since the image of $j$ is contained in $E_0^*$, Lemma~\ref{lem.tue-2}
tells us that
$\Alg(j)$ is defined over $E_0^{\Sym_n}$. But
$E_0 \simeq k(P)$
and so $E_0^{\Sym_n} \simeq k(P)^{\Sym_n}$ which is indeed rational over
$k$. This completes the proof of Theorem~\ref{thm.rat}(a).
\qed

\subsection{Proof of part(b)} By the Merkurjev-Suslin Theorem,
\[ \Mat_r(A) = (a_1, b_1)_{n_1} \otimes_F \dots
\otimes_F (a_l, b_l)_{n_l} \, , \]
for some $r, l \geq 1$, where
$(a, b)_n$
denotes a symbol algebra; see \eqref{E:symbol}.

Let $\lambda_1, \dots, \lambda_l, \mu_1, \dots, \mu_l$ be $2l$
central variables, algebraically independent over $F$.
We will write $\lambda$ in place of $(\lambda_1, ..., \lambda_l)$
and $\mu$ in place of $(\mu_1, ..., \mu_l)$. Then
\begin{eqnarray*}
\Mat_r(A)(\lambda, \mu) =
(a_1, b_1)_{n_1} \otimes_{K(\lambda, \mu)} \dots
\otimes_{K(\lambda, \mu)}(a_l, b_l)_{n_l} =  \\
(a_1', b_1')_{n_1} \otimes_{K(\lambda, \mu)} \dots
\otimes_{K(\lambda, \mu)} (a_l', b_l')_{n_l} = \\
\Bigl( (a_1', b_1')_{n_1} \otimes_{F_0} \dots
\otimes_{F_0} (a_l', b_l')_{n_l} \Bigr) \otimes_{F_0} K(\lambda, \mu)
\, ,
\end{eqnarray*}
where $a_i' = a_i \lambda_i^{n_i}$ and $b_i' = b_i \mu_i^{n_i}$ for
$i = 1, \dots, l$ and
$F_0 = k(a_1', b_1', ..., a_l', b_l')$.
This shows that $\Mat_r(D)(\lambda, \mu)$ is defined over $F_0$.
It remains to prove that $F_0$ is rational over $k$.  The $2l$ elements
$a_1', b_1', \dots, a_l', b_l'$ are clearly algebraically independent over
$F$. Hence, they are algebraically independent over $k$, and consequently,
$F_0$ is rational over $k$, as claimed.
\qed

\begin{remark} \label{rem.explicit-bound}
Our proof of Theorem~\ref{thm.rat} can be used to deduce explicit
lower bounds on $s$ in parts (a) and (b) from explicit lower bounds
in Theorems of Bessenrodt-LeBruyn~\cite[Proposition 3]{bl} (on $\rank(L)$)
and Merkurjev-Suslin (on $r$). The lowest possible value of $r$ in
part (b), called the {\em Merkurjev-Suslin number}, is of independent
interest; see~\cite[Section 7.2]{rowen-rt2}.
\end{remark}

\section{Proof of Theorem~\ref{thm.ab-cr}}
\label{sect.deg4}

\begin{reduction} \label{red6.1}  In the course of proving
Theorem~\ref{thm.ab-cr}, we may assume without loss of generality
that $A$ is a division algebra.
\end{reduction}

Indeed, let $D = \Alg(\mu)$, as in Theorem~\ref{thm.tue-3}, with
$\cG = \bbZ/m \times \bbZ/2$, and $\cH = \{ 1 \}$. Then $D$ is
a division algebra (see Remark~\ref{rem.exponent}), and any other
$\cG$-crossed product $A/F$ is a rational specialization of $D$.
Thus, if we know that Theorem~\ref{thm.ab-cr} holds for $D$ then
it holds for $A(t_1, \dots, t_s)$, where $t_1, \dots, t_s$
are independent variables over $F$.
Using induction on $s$,
we see that Reduction~\ref{red6.1} is now a consequence
of the following lemma (applied to $B = \Mat_m(A)$, with $r = 2$,
$m_1 = m$ and $m_2 = 2m$):

\begin{lem} \label{lem6.05}
Let $B/K$ be a central simple algebra of degree $d = m_1 \dots m_r$
and let $t$ be an independent variable over $K$. Assume $K$ contains
a primitive root of unity of degree lcm$(m_1, \dots, m_r)$. If
\[ B(t) = (a_1(t), b_1(t))_{m_1} \otimes \dots \otimes
(a_r(t), b_r(t))_{m_r} \]
for some $a_1(t), b_1(t), \dots, a_r(t), b_r(t) \in K(t)$ then
\[ B = (a'_1, b'_1)_{m_1} \otimes \dots \otimes (a'_r, b'_r)_{m_r} \]
for some $a'_1, b'_1, \dots, a'_r, b'_r \in K$.
\end{lem}

Our proof is based on a standard specialization argument;
for the sake of completeness, we supply the details below.

\begin{proof} We may assume that $K$ is an infinite field; otherwise
$B$ is a matrix algebra over $K$, and we can take, e.g., $a'_i = 1$,
$b'_i = -1$ for every $i$.

Choose generators $x_i(t)$ and $y_i(t)$ for the cyclic subalgebra
$(a_i(t), b_i(t))_{m_i}$ of $B(t) = B \otimes_K K(t)$ such that
$x_i(t)^{m_i} = a_i(t)$, $y_i(t)^{m_i} = b_i(t)$, and
$x_i(t) y_i(t) = \zeta_{m_i} y_i(t) x_i(t)$, where $\zeta_{m_i}$
is a primitive root of unity of degree $m_i$ in $K$.
Choose a $K$-basis $b_1, \dots, b_{d^2}$ of $B$ and write
\begin{equation} \label{e6.05}
x_i(t) = \sum_{j = 1}^{d^2} \alpha_{ij}(t) b_i \quad \text{and}
\quad
y_i(t) = \sum_{j = 1}^{d^2} \beta_{ij}(t) b_i \, , \end{equation}
for some $\alpha_{ij}(t), \beta_{ij}(t) \in K(t)$.
Since $K$ is an infinite field, we can choose $t_0 \in K$ such that
$\alpha_{ij}(t_0)$ and $\beta_{ij}(t_0)$ are well-defined and
\[ x_i(t_0) = \sum_{j = 1}^{d^2} \alpha_{ij}(t_0) b_i \quad \text{and}
\quad
y_i(t_0) = \sum_{j = 1}^{d^2} \beta_{ij}(t_0) b_i \]
are non-zero.  Let $B_i$ denote the subalgebra of $B$ that is
generated by $x_i(t_0)$ and $y_i(t_0)$.
Then $B_i = (a'_i, b'_i)_{m_i}$, where
$a'_i = x_i(t_0)^{m_i} = a_i(t_0)$ and $b'_i = y_i(t_0)^{m_i} = b_i(t_0)$,
and $B_1, \dots, B_r$ are commuting subalgebras of $B$ of degrees
$m_1, \dots, m_r$. Hence, by the double centralizer theorem
(cf., e.g., \cite[Theorem 12.7]{pierce}), $B = B_1 \otimes \dots \otimes
B_r$.
This completes the proof of Lemma~\ref{lem6.05} and thus
of Reduction~\ref{red6.1}.
\end{proof}

We now continue with the proof of Theorem~\ref{thm.ab-cr}.  In the course
of the proof we shall use the following notations. Write
$\cG = \bbZ/m \times \bbZ/2 = \langle\sigma_1, \sigma_2\rangle$,
where $\sigma_1^{m} = \sigma_2^2 = 1$. Let $K = F(\alpha_1, \alpha_2)$,
be a maximal $\cG$-Galois subfield of $A$, where $\alpha_1^m = a_1$
and $\alpha_2^2 = a_2$ are elements of $F$, and
\begin{align}
\sigma_1(\alpha_1) &= \zeta_m \alpha_1 \, , &
\sigma_1(\alpha_2) &= \alpha_2 \, ,  & \nonumber \\
\sigma_2(\alpha_1) &= \alpha_1 \, , &
\sigma_2(\alpha_2) &= - \alpha_2 \, .
\label{e.alpha}
\end{align}
Here $\zeta_m \in F$ is a primitive $m^{\text{th}}$ root of unity, so that
$K$ is, indeed, a $\cG$-Galois extension of $F$. Note that
the statement of Theorem~\ref{thm.ab-cr} assumes that $F$ contains
not only a primitive $m$th root of unity $\zeta_m$ but also
a primitive $2m^{\text{th}}$ root of unity $\zeta_{2m}$; we shall make
use of $\zeta_{2m}$ later in the proof.

By the Skolem-Noether Theorem,
there exist units $z_1, z_2 \in A$ such that $z_i x z_i^{-1} = \sigma_i(x)$
for every $x \in K$ ($i = 1,2$). Set
\begin{equation} \label{e.z1-z2}
z_1^m = b_1 \in F(\alpha_2)^*,\ z_2^2 = b_2 \in F(\alpha_1)^*,\
\text{and} \ u = z_1z_2z_1^{-1}z_2^{-1} \in K^*\ .
\end{equation}
By~\cite[Theorem 1.3]{as},
the algebra structure of $A$ can be recovered from the $\cG$-field $K$
and the elements $u \in K^*$, $b_1 \in F(\alpha_2)^*$ and
$b_2 \in F(\alpha_1)^*$. (These elements have to satisfy certain
compatibility conditions; the exact form of these
conditions shall not concern us in the sequel.)
We will write $A = (K, \cG, u, b_1, b_2)$.

\begin{lem} \label{lem.tensor}
Let $A = (K, \cG, u, b_1, b_2)$ and $A' = (K, \cG, u', b_1', b_2')$
be $\cG$-crossed products.  Then $A \otimes_F A'$ is Brauer equivalent
to $(K, \cG, uu', b_1 b_1' ,b_2 b_2')$.
\end{lem}

\begin{proof}
The class of $A = (K, \cG, u, b_1, b_2)$ in the relative Brauer group
$\BB(K/F) = H^2(\cG,K^*)$
is given by a normalized $2$-cocycle $a \colon \cG \times \cG \to K^*$ so
that
$$
b_i = a(\sigma_i,\sigma_i)a(\sigma_i^2,\sigma_i)\dots
a(\sigma_i^{m_i-1},\sigma_i)
$$
holds for $i = 1,2$, where $m_1 = m$ and $m_2 = 2$, and
$$
u = a(\sigma_1,\sigma_2)a(\sigma_s,\sigma_1)^{-1}\ .
$$
Similarly, the class of $A'$ is given by a $2$-cocycle $a'$. Then the class
of $A\otimes_F A'$ is given by the cocycle $a a'$;
see, e.g., \cite[Proposition 14.3]{pierce}. This proves the lemma.

The following alternative ring-theoretic argument was suggested
by the referee:
Choose $z_1, z_2 \in A$, as in~\eqref{e.z1-z2}, and similarly
for $z_1'$, $z_2'$ in $A'$.
The subalgebra $S$ of $A \otimes_F A'$ generated by $K \otimes 1 $,
$z_1 \otimes z_1'$ and $z_2 \otimes z_2'$, is clearly isomorphic to
$(K, \cG, u u', b_1 b_1', b_2 b_2')$. Its centralizer $C_{A \otimes A'}(S)$
is
an $F$-central simple algebra of degree $2m = [K:F]$, containing
$(K \otimes K)^{\cG}$, where $\cG$ acts diagonally
on $K \otimes K$. Since $\cG$ is abelian,
$(K \otimes K)^{\cG} \simeq F \oplus \dots \oplus F$, as an
$F[\cG]$-algebra.
In particular, $C_{A \otimes A'}$ contains the idempotents of $K \otimes K$
and, hence, is split over $F$. We thus conclude that
\[ A \otimes_F A' \simeq S \otimes_F C_{A \otimes_F A'} (S) \sim S \simeq
(K, \cG, u u', b_1 b_1' ,b_2 b_2') \, , \]
as claimed. (Here $\sim$ denotes Brauer equivalence over $F$.)
\end{proof}

\smallskip
We now proceed with the proof of Theorem~\ref{thm.ab-cr}, using
the notations of~\eqref{e.alpha} and~\eqref{e.z1-z2}.
Since $b_1 = z_1^m \in K^{\sigma_1} = F(\alpha_2)$, we
can write
\begin{equation} \label{e.b1}
b_1 = f_1 + f_2 \alpha_2 \, ,
\end{equation}
for some $f_1, f_2 \in F$.

\begin{lem} \label{lem.ab-cr}
\begin{enumerate}

\item
If $f_1 = 0$ then $A$ is cyclic.
\item
If $f_2 = 0$ then $A = (a, b)_m \otimes (c, d)_2$, for some $a, b, c, d
\in F^*$.
\end{enumerate}
\end{lem}

\begin{proof}
(a) If $f_1 = 0$ then $z_1^{2m} = b_1^2 = f_2^2 a_2 \in F^*$ but
$z_1^m = f_2 \alpha_2 \not \in F$.  Since $F$ contains
a primitive root of unity of degree $2m$, $F(z_1)$ is a cyclic maximal
subfield of $A$ of degree $2m$; cf.~\cite[Theorem VIII.6.10(b)]{lang}.
Thus $A$ is a cyclic algebra, as claimed.

(b) If $f_2 = 0$, i.e., $b_1 \in F$, then the $F$-subalgebra $A_0$
of $A$ generated by $z_1$ and $\alpha_1$ is cyclic of degree $m$.
By the double centralizer theorem, $A = A_m \otimes Q$,
where $Q$ is a quaternion algebra, as claimed.
\end{proof}

We are now ready to finish the proof of Theorem~\ref{thm.ab-cr}.
Lemma~\ref{lem.ab-cr} tells us that Theorem~\ref{thm.ab-cr} is immediate
if $f_1 = 0$ or $f_2 = 0$. Thus from now on we shall assume $f_1 f_2 \neq
0$.

Now let $A = (K, \cG, u, b_1, b_2)$ and, for any  $f \in F^*$, define
$A_f = (K, \cG, u, f b_1, b_2)$.
Since $(a_1, f)_m \otimes_F \Mat_2(F) \simeq (K, \cG, 1, f, 1)$,
Lemma~\ref{lem.tensor} tells us that $A_f \sim
(a_1, f)_m \otimes_F A$, where $\sim$
denotes Brauer equivalence. In other words,
$A \sim (f, a_1)_m \otimes_F A_f$.
Thus it is enough to show that
$A_f$ is cyclic, for some $f \in F^*$.

To prove the last assertion, observe that
if we expand $(z_1 + \alpha_1)^m$ then all terms,
other than $z_1^m$ and $\alpha_1^m$, will cancel. (For a simple proof
of this fact, due to Bergman, see~\cite[p.~195]{rowen-rt2}). Thus,
if $\gamma = z_1 + \alpha_1$ then
\[ \gamma^m = z_1^m + \alpha_1^m = f b_1 + a_1 \, \]
in $A_f$.  Setting $f = -\frac{a_1}{f_1} \in F^*$, we obtain
$\gamma^m = c \alpha_2$, where $c = -\frac{a_1 f_2}{f_1} \in F^*$;
see~\eqref{e.b1}.
Thus $\gamma^{2m} = c^2 a_2 \in F^*$ but $\gamma^m \not \in F$.
Since $F$ contains a primitive $2m^{\text{th}}$ root of unity,
$F(\gamma)/F$ is a cyclic field extension
of degree $2m$. In other words, $F(\gamma)$ a cyclic maximal
subfield of $A_f$, and $A_f$ is a cyclic algebra of degree $2m$, as claimed.
\qed

\section{The field of definition of a quadratic form}
\label{sect.quadr}

\subsection{Preliminaries}
Let $V = F^n$ be an $F$-vector space, equipped with a quadratic
form $q \colon V \to F$. Recall that $q$ is said to be defined over
a subfield $F_0$ of $F$ if $q = q_{F_0} \otimes F$, where $q_{F_0}$
is a quadratic form on $V_0 = F_0^n$.  Is easy to see that $q$
is defined over $F_0$ if and only if $V$ has an $F$-basis
$e_1, \dots, e_n$ such that $b(e_i, e_j) \in F_0$, where
$b \colon V \times V \to F$ is the symmetric bilinear form
associated to $q$ (i.e., $q(v) = b(v, v)$).

We shall always assume that $\Char(F) \neq 2$ and
$F$ (and $F_0$) contain a base subfield $k$.
As usual, we shall write $\mo a_1, \dots, a_n \mc$ for the diagonal form
\[ a_nx_1^2 + \dots + a_n x_n^2 \]
and $\ll a_1, \dots, a_n \gg$ for the Pfister
form $\mo 1, a_1 \mc \otimes \dots \otimes \mo 1, a_n \mc$.
Given a quadratic form $q$ we shall ask:

\smallskip
(a) What is the smallest value of $\trdeg_k(F_0)$, where $q$ is defined
over $F_0$? We shall denote this number by $\tau(q)$.

\smallskip
(b) Can $q$ be defined over a rational extension $F_0$ of $k$?

\smallskip
\noindent
These are the same questions we asked for central simple algebras in
the Introduction.  In the case of quadratic forms our answers are more
complete (and the proofs are easier).

\begin{prop} \label{prop7.1} Let $V = F^n$ and let $q \colon V
\to F$ be a quadratic form on $V$. Then
\begin{enumerate}
\item
$\tau(q) \leq n$.
Moreover, if $a_1, \dots, a_n$ are independent variables over $k$,
$F = k(a_1, \dots, a_n)$, and $q = \mo a_1, \dots, a_n \mc$ then
$\tau(q) = n$.
\item
Let $t_1, \dots, t_n$ be independent variables over $F$.
Then $q' = q \otimes_F F(t_1, \dots, t_n)$ is defined over
a rational extension $F_0$ of $k$.
\end{enumerate}
\end{prop}

\begin{proof} Diagonalizing $q$, write $q = \mo a_1, \dots, a_n \mc$
in the basis $e_1, \dots, e_n$.

\smallskip
(a) To prove the first assertion, set $F_0 = k(a_1, \dots, a_n)$.
Then $q$ is defined over $F_0 = k(a_1, \dots, a_n)$ and
$\trdeg_k (F_0) \leq n$, as desired. For the proof of the second
assertion see~\cite[Proof of Theorem 10.3]{reichstein1}.

\smallskip
(b) Set $a_i' = t_i^2 a_i$. Then $q' = \mo a_1, \dots, a_n \mc =
\mo a_1', \dots, a_n' \mc$ over $F(t_1, \dots, t_n)$. Hence,
$q'$ is defined over $F_0 = k(a_1', \dots, a_n')$.
We claim that $F_0$ is rational over $k$. Indeed, since the nonzero
elements of $\{ a_1', \dots, a_n' \}$ are algebraically
independent over $F$, they are algebraically independent
over $k$, and the claim follows.
\end{proof}

In the sequel we shall need the following analogue of Lemma~\ref{rattau}.

\begin{lem} \label{rattau1} Let $q$
be a quadratic form defined over $F$, $t_1, \dots, t_r$ be independent
variables over $F$, and $F' = F(t_1, \dots, t_r)$.
Set $q' = q \otimes _F F'$.  Then $\tau(q) = \tau(q')$.
\end{lem}

\begin{proof} The inequality $\tau(q') \leq \tau(q)$ is obvious from
the definition of $\tau(q)$. To prove the opposite inequality,
we may assume $F$ is an infinite field; otherwise $\tau(q) =0$,
and there is nothing to prove. We may also assume
$r = 1$; the general case
will then follows by induction on $r$.
Let $b'$ be the symmetric bilinear form associated to $q'$
and choose a basis $b_1(t), \dots, b_n(t)$ of $(F')^n$ so
that $\trdeg_k \, k( \alpha_{ij}(t)) = \tau(q')$, where $\alpha_{ij}(t)
= b'(b_i(t),b_j(t))$.
Since $F$ is an infinite field, we can find a $c \in F$ such that
(i) the vectors $b_1(c), \dots, b_d(c)$ are
well-defined and form a basis of $F^d$, and
(ii) each $\alpha_{ij}(c)$ is well-defined. Now
$q$ is defined over $k(\alpha_{ij}(c))$ and thus
\[ \tau(q) \leq \trdeg_k(\alpha_{ij}(c)) \leq
\trdeg_k(\alpha_{ij}(t)) = \tau(q') \, , \]
as claimed.
\end{proof}

\subsection{Proof of Theorem~\ref{thm.trace-form}}

Let $F$ be a field containing a primitive $4^{\text{th}}$ root
of unity. Note that for the purpose of proving
Theorem~\ref{thm.trace-form}, we may assume that $A/F$ is a division
algebra. Otherwise, $A$ is isomorphic to $\Mat_4(F)$ or to
$\Mat_2(D)$, where
$D = (a, b)_2$ is a quaternion algebra. Thus, $A$ is defined over
$k$ or over the field $F_0 = k(a, b)$, respectively, and so is the trace
form
of $A$.
Alternatively, a simple direct computation shows that the trace
form of $\Mat_2(E)$ is trivial (and thus is defined over $k$)
for any central simple algebra $E/F$.

From now on we will assume that $A/F$ is a division algebra of degree 4.
By a theorem of Albert~\cite{albert}, $A$
is a $\cG$-crossed product, with $\cG = \bbZ/2 \times \bbZ/2$.
Let $K$ be a $\cG$-Galois
maximal subfield. Using the notations introduced in
Section~\ref{sect.deg4} (with
$m = 2$), we write $\cG = \mo \sigma_1, \sigma_2 \mc$, $K = F(\alpha_1,
\alpha_2)$, $\alpha_i^2 = a_i \in F$ and
$A = (K, \cG, u, b_1, b_2)$ for some $u \in K^*$,
$b_1 \in F(\alpha_2) = K^{\sigma_1}$, and
$b_2 \in F(\alpha_1) = K^{\sigma_2}$.
Set $\sigma_3 = \sigma_1 \sigma_2 \in \Gal(K/F)$; $z_3 = (z_1z_2)^{-1}$,
$\alpha_3 = \alpha_1 \alpha_2$, $a_3 = \alpha_3^2 = a_1 a_2$,
$b_3 = z_3^2$ (so that $b_i = z_i^2$ for $i = 1, 2, 3$), and
\[ t_i = \frac{1}{2}\tr_{K^{\sigma_i}/F}(z_i^2) \, , \quad
n_i = N_{K^{\sigma_i}/F}(z_i^2) \,  \]
for $i = 1, 2, 3$.

Our proof of Theorem~\ref{thm.trace-form}
is based on the following result of Serre~\cite{serre},~\cite{rst}.

\begin{prop} \label{prop.serre}
Suppose $z_1$ and $z_2$ are chosen so that $t_i \neq 0$
and $n_i^2 - t_i \neq 0$ for any $i = 1, 2, 3$. Then
the trace form $q$ of $A$ is Witt-equivalent (over $F$)
to $q_2 \oplus q_4$, where
\[ q_2 = \ll n_1 - t_1^2, \; n_2 \gg \]
is a 2-fold Pfister form and
\[ q_4 = \ll t_1 - n_1^2, \;
(n_2 - t_2^2) n_2, \; t_1 t_2, \;  t_2 t_3 \gg \]
is a 4-fold Pfister form.
\qed
\end{prop}

We claim that for the purpose of proving Theorem~\ref{thm.trace-form}, we
may assume without loss of generality that $t_i \neq 0$
and $n_i^2 - t_i \neq 0$ for any $i = 1, 2, 3$.
Indeed, in view of Lemma~\ref{rattau1} it suffices to
prove~Theorem~\ref{thm.trace-form} for a single division
algebra $A$ which has the rational specialization property
in the class of algebras of degree 4, e.g., for $A = \UD(4)$;
see Remark~\ref{rem-ud(n)}. Thus we only need to show that
in this algebra $t_i \neq 0$ and $n_i - t_i^2 \neq 0$
for any choice of $z_1, z_2$.

Indeed, we may assume without loss of generality that $i = 1$
(the cases where $i = 2$ and $3$ will then follow by symmetry).
Write $b_1 = f_1 + f_2\alpha_2$ for some $f_1, f_2 \in F$, where
$t_1 = f_1$ and $n_1 - t_1^2 = f_2^2 a_2$.
Lemma~\ref{lem.ab-cr} shows that if $t_1 = 0$
then $A$ is cyclic and if $n_1 - t_1^2 = 0$ then $A$ is biquaternion.
But since our algebra $A$ has the rational specialization property,
it is neither cyclic nor biquaternion. We
conclude that $t_1(n_1 - t_1^2) \neq 0$, as claimed.

We now proceed to simplify the form given by Proposition~\ref{prop.serre}.
After expanding $q_2$ and $q_4$, cancelling the common
term $\mo 1, t_1^2 - n_1 \mc$ (which can be done, since
we are assuming $\sqrt{-1} \in F$) and dividing some of the entries
by elements of $(F^*)^2$,
we see that the trace form of $A$ is Witt equivalent to
the 16-dimensional form
\begin{equation} \label{e.equiv-form}
q = \mo 1 \, , \,  1 - \frac{n_1}{t_1^2} \mc \, \otimes
(\mo \frac{n_2}{t_2^2}
\mc \oplus
\pfistero (1 - \frac{n_2}{t_2^2})\frac{n_2}{t_2^2} \, , \;
t_1 t_2, \; t_2 t_3 \pfisterc _0 )\, ,
\end{equation}
where $\pfistero \lambda_1, \, \, \dots, \lambda_r \pfisterc _0$
is defined as a $2^r-1$-dimensional form such that \\
$\pfistero \lambda_1, \, \lambda_2, \, \dots, \lambda_r
\pfisterc _0 \oplus \mo 1 \mc$
is the $n$-fold Pfister form
$\pfistero \lambda_1, \, \lambda_2, \, \dots, \lambda_r \pfisterc$.

Note that since $q$ and the trace form of $A$
are Witt equivalent 16-dimensional forms,
the Witt Decomposition Theorem implies that they are, in fact, the same
(i.e., isometric).  We now observe that all entries of $q$
lie in the subfield $F_0 = k(\frac{n_1}{t_1^2},
\frac{n_2}{t_2^2}, t_1 t_2, t_2 t_3)$ of $F$. Thus
the trace form of $A$ is defined over $F_0$.
Since $F_0$ is generated by 4 elements over $k$,
we have $\trdeg_k \, F_0 \leq 4$. This completes the proof
of Theorem~\ref{thm.trace-form}.
\qed

\end{document}